\documentclass[a4paper,reqno,10pt,oneside]{amsart}
\usepackage{amsmath}
 \usepackage[
  hmarginratio={1:1},     % equal left and right margins
  vmarginratio={1:1},     % equal top and bottom margins
  textwidth=14cm,        % new text width  
  textheight=22cm,  
  heightrounded,          % always useful
]{geometry}

\usepackage{graphicx}
\usepackage{mathrsfs}
\usepackage{amssymb,amsmath, amsfonts, amsthm}
\usepackage{latexsym}
\usepackage{color} %{\color{red}..............OK}
\usepackage[dvips]{epsfig}
\usepackage{graphicx}

\allowdisplaybreaks[1]

\numberwithin{equation}{section}

\renewcommand{\(}{\left(}
\renewcommand{\)}{\right)}
\renewcommand{\[}{\left[}
\renewcommand{\]}{\right]}

\newlength{\defbaselineskip}
\newtheorem{theorem}{Theorem}[section]
\newtheorem{proposition}[theorem]{Proposition}
\newtheorem{lemma}[theorem]{Lemma}

\renewcommand{\le}{\leqslant}
\renewcommand{\ge}{\geqslant}

\newcommand{\beq}{\begin{equation}}
\newcommand{\eeq}{\end{equation}}
\newcommand{\beqs}{\begin{equation*}}
\newcommand{\eeqs}{\end{equation*}}
\newcommand{\beqn}{\begin{eqnarray}}
\newcommand{\eeqn}{\end{eqnarray}}
\newcommand{\beqns}{\begin{eqnarray*}}
\newcommand{\eeqns}{\end{eqnarray*}}
\newcommand{\bdoc}{\begin{document}}
\newcommand{\edoc}{\end{document}}
\newcommand{\be}{\begin{enumerate}}
\newcommand{\ee}{\end{enumerate}}
\newcommand{\bdescr}{\begin{description}}
\newcommand{\edescr}{\end{description}}
\newcommand{\ba}{\begin{array}}
\newcommand{\ea}{\end{array}}
\newcommand{\intR}{\int_{\mathbb R^N}}
\newcommand{\R}{\mathbb R^N}
\newcommand{\e}{\varepsilon}
 \renewcommand{\(}{\left(}
\renewcommand{\)}{\right)}
\renewcommand{\[}{\left[}
\renewcommand{\]}{\right]}
\newenvironment{Proof}{\noindent{\bf Proof}}{\hfill$\Box$\\[2mm]}

\def\red{\color{red}}
\def\bl{\color{blue}}
\def\bk{\color{black}}
\begin{document}

\title[Nondegeneracy of the bubble for the  critical $p$-Laplace equation]{Nondegeneracy of the bubble for the critical $p-$Laplace equation}

\author{Angela Pistoia and Giusi Vaira}
\address{Angela Pistoia\\ Dipartimento di Scienze di Base e Applicate per l'Ingegneria, Sapienza Universit\`a di Roma\\Via Scarpa 16, 00161 Roma, Italy}
\email{angela.pistoia@uniroma1.it}
\address{Giusi Vaira\\ Dipartimento di Matematica e Fisica, Universit\`a degli studi della Campania ``Luigi Vanvitelli"\\Viale Lincoln 5, 81100 Caserta, Italy}
\email{giusi.vaira@unicampania.it}
\maketitle
\begin{abstract}
We prove the non-degeneracy of the extremals of the 
Sobolev inequality 
$$\int\limits_{\mathbb R^N}|\nabla u|^pdx\ge \mathcal S_p\int\limits_{\mathbb R^N}|u|^{Np\over N-p}dx,\ u\in \mathcal D^{1,p}(\mathbb R^N)$$
when $1<p<N,$ as solutions of a critical  quasilinear equation involving the $p-$Laplacian.

\end{abstract}
\section{Introduction and statement of the main result}
In this paper we establish the {\it linear non-degeneracy} of the extremals of the optimal classical Sobolev inequality
\begin{equation}\label{sobolev}
\mathcal S _p \|u\|_{L^{p^*} (\mathbb R^N)}\le \|\nabla u\|_{L^p(\mathbb R^N)}\ \hbox{for any}\ u\in\mathcal D^{1,p}(\mathbb R^N),
\end{equation}
where $p^* :=\frac{Np}{N-p}$ and $1<p<N$.\\

Aubin \cite{aubin} and Talenti \cite{Talenti}  found the optimal constant and the extremals for inequality \eqref{sobolev}. Indeed, equality is achieved precisely by the functions
 \begin{equation}\label{solcrit}
 U_{\delta, \xi}(x):=\delta^{-\frac{N-p}{p}} U\left(\frac{x-\xi}{\delta}\right)\  \hbox{where}\ \delta>0,\ \xi\in\mathbb R^N
\end{equation}
where
 \begin{equation}\label{solcrit1}
 U(x)=\left(\frac{\alpha_{N, p}}{1+|x-\xi|^{\frac{p}{p-1}}}\right)^{\frac{N-p}{p}}\  \hbox{with}\ \alpha_{N, p}:=N^{\frac 1 p}\left(\frac{N-p}{p-1}\right)^{\frac{p-1}{p}},
\end{equation}
which solve the critical equation
\begin{equation}\label{equ}
-\Delta _p u=u^{p^* -1}\ \hbox{in}\ \mathbb R^N,\ u>0\ \hbox{in}\ \mathbb R^N,\  u\in\mathcal D^{1,p}(\mathbb R^N).
\end{equation}

All the solutions to the equation \eqref{equ} are indeed the only ones of \eqref{solcrit}. Caffarelli, Gidas and Spruck proved the claim when $p=2$. The case $p\not=2$ has been firstly solved by Guedda and Veron \cite{gv1} in the radial case, where the authors classified all the positive radial solutions  and successively by Damascelli, Merch\'an, Montoro and Sciunzi \cite{DMMS} when ${2N\over N+2}\le p<2$, by Vet\'ois \cite{Vetois} and  Damascelli and Ramaswamy \cite{DR} when $1<p<2$ and finally by Sciunzi \cite{S} in the remaining cases, namely when $2<p<N$.\\

Here we are interested in the linear non-degeneracy of the solutions \eqref{solcrit} to equation \eqref{equ}.

Let us point out that equation \eqref{equ} is invariant by scaling and by translations. Therefore, if we differentiate the equation
$$-\Delta_p U_{\delta,\xi}= U_{\delta,\xi}^{p^*-1}\ \hbox{in}\ \mathbb R^N$$
with respect to the parameters $\delta $ and $\xi_1,\dots,\xi_N$ at $\delta=1$ and $\xi=0$ we see that the functions
 \begin{equation}\label{Z0}Z_0(x):=-\partial_\delta  U_{\delta, \xi}|_{\delta=1, \xi=0}=\frac{N-p}{p} U+x\cdot \nabla  U \end{equation} and \begin{equation}\label{Zi}Z_i(x):=\partial_{\xi_i}  U_{\delta, \xi}|_{\delta=1, \xi=0}=-\partial_{x_i} U,\  i=1, \ldots, N\end{equation} annihilate the linearized operator around the function $ U$ defined in \eqref{solcrit1}, namely they solve the linear equation
\begin{equation}\label{linpb}
  -{\rm{div}}\(|\nabla  U|^{p-2}\nabla \phi\)-(p-2){\rm{div}}\(|\nabla  U|^{p-4}\(\nabla  U,\nabla\phi\)\nabla  U\)=(p^*-1) U^{p^*-2}\phi\ \hbox{in}\ \mathbb R^N.
  \end{equation}

  We say that $ U $ is non-degenerate if the kernel of the associated linearized operator \eqref{linpb} is spanned only by the  functions $Z_i$'s defined in \eqref{Z0} and \eqref{Zi}. This property is true when $p=2$ as it was established by Rey in \cite{r}.
  Our main result extends the non-degeneracy of the solution $ U$  to any $p\in (1,N)$ in    the weighted Sobolev space     $\mathcal D^{1, 2}_*(\mathbb R^N)$, which is defined as the completion of $C^1_c(\mathbb R^N)$ with respect to the norm
$$\|\phi\|:=\left( \int_{\mathbb R^N}|\nabla U|^{p-2}|\nabla\phi|^2\, dx\right)^{1/2}.$$
   \begin{theorem}\label{main1}
The solution 
$$ U(x)=\left(\frac{\alpha_{N, p}}{1+|x |^{\frac{p}{p-1}}}\right)^{\frac{N-p}{p}}\  \hbox{with}\ \alpha_{N, p}:=N^{\frac 1 p}\left(\frac{N-p}{p-1}\right)^{\frac{p-1}{p}}$$
 of equation \eqref{equ} is non-degenerate in the sense that all the solutions  of the equation \eqref{linpb} in the space $\mathcal D^{1, 2}_*(\mathbb R^N)$
  are linear combination of the functions 
 $$Z_0(x)=\frac{N-p}{p} U+x\cdot \nabla  U,\ Z_1(x)=\partial_{x_1}  U(x),\dots,\ Z_N(x)=\partial_{x_N}  U(x).$$
 \end{theorem}
%
%Actually, a crucial ingredient in studying bubbling concentration behaviour of families of solutions to  elliptic equations involving critical Sobolev growth is understading the linear non-degeneracy of solutions \eqref{solcrit} to problem \eqref{equ}. 

\medskip
 The structure of the linearized equation \eqref{linpb} strongly suggests to introduce the   space  $\mathcal D^{1, 2}_*(\mathbb R^N)$. A similar first order Sobolev space with weight was introduced by Damascelli and Sciunzi in \cite{DS} to study a linearized operator on a bounded domain. Here the situation is much more delicate due to the unboundness of the domain. Section \ref{sec-1}   is devoted to prove some   properties of $\mathcal D^{1, 2}_*(\mathbb R^N)$ which are essential to get Theorem \ref{main1} whose proof is carried out in Section \ref{sec-2}. \\

Quasilinear equations with critical growth involving the $p-$Laplace operator have been widely studied in recent years using a variational framework, starting from the quasilinear version of the classical Brezis-Nirenberg problem (see \cite{bn})   studied by Guedda and Veron in \cite{gv}.
 In particular, we would like to focus on  the problem of the existence of sign-changing solutions to the critical equation
\begin{equation}\label{cri}-\Delta _p u=|u|^{p^*-2}u\ \hbox{in}\ \Omega,
\end{equation}
where  $\Omega$ is  either the whole space $\mathbb R^N$ or     a bounded smooth domain in $\mathbb R^N$ in which case we assume homogeneous Dirichlet boundary conditions.
As far as we know the only result concerning existence of sign-changing solutions to \eqref{cri} in the whole space is due to Clapp and Lopez Rios in \cite{cl},
where they prove that \eqref{cri} has a certain finite number (depending on the dimension $N$) of non-radial sign-changing solutions.
On the other hand if $p=2$  del Pino, Musso, Pacard, and Pistoia in \cite{dmpp1,dmpp2}  used the Lyapunov–Schmidt
procedure to build infinitely many sign-changing solutions which look like a positive bubble crowned by an arbitrary large number of negative   bubbles arranged on a regular polygon.  It would be interesting to check if it is possible to build this kind of solutions in the quasilinear case.
When $\Omega$ is a bounded domain, the existence of solutions is a more delicate issue.
 Indeed if $\Omega$ is starshaped the problem does not have any solutions because of a Pohozaev identity obtained by Guedda and Veron in \cite{gv}. 
The existence of a positive solution  has been proved  by Mercuri, Sciunzi and Squassina in \cite{mss} when the domain has a small hole, in the same spirit of Coron's result  \cite{c} when $p=2.$  The existence of a sign-changing solution has been
   obtained by Mercuri and Pacella in \cite{mp} when the domain $\Omega$  has either a small hole and little symmetry  or    a hole of any size and more symmetry.   On the other hand, if $p=2$ and $\Omega$ has a small hole,      Musso and Pistoia in \cite{mup} (see also \cite{gmp,gmpp})  used the Lyapunov–Schmidt
procedure to  built sign-changing solutions which look like   the superposition  of   bubbles with alternating sign whose number
  becomes arbitrary large as the size of the hole approaches zero. It is natural to ask if this kind of solutions do exist also in the quasilinear case.
  
In both cases the understanding of the linear non-degeneracy of the bubble is the first step in the application of the Ljapunov-Schmidt procedure.\\

{\bf Acknowledgements.} We wish to thank professor Berardino Sciunzi for many helpful comments and discussions.
Moreover, we warmly thank the anonymous referee for his/her valuable comments which allow us to improve the presentation of the paper.\\

\section{A suitable weighted Sobolev space}\label{sec-1}
 
First of all, let us point out the   following fact.
\begin{lemma}\label{mah}
\begin{itemize}
\item[(i)]
If $p\in(1,2)$ there exists $C>0$ such that  
$$\(\int\limits_{\mathbb R^N} |\phi|^{Np\over N-p}\, dx\)^{N-p\over p}\le C \(\int\limits_{\mathbb R^N} |\nabla U|^{p-2}|\nabla \phi|^{2}\, dx\)^{\frac12} \hbox{for any}\ \phi\in C^1_c(\mathbb R^N).$$
\item[(ii)]If $p\in(2,N)$   for any $R>0$ there exists $C(R)>0$ such that
$$\(\int\limits_{\mathbb R^N\setminus B_R(0)}|x|^{-{Np+p-2N\over p-1}} |\phi|^2\, dx\)^{\frac12}\le C(R) \(\int\limits_{\mathbb R^N} |\nabla U|^{p-2}|\nabla \phi|^{2}\, dx\)^{\frac12} \hbox{for any}\ \phi\in C^1_c(\mathbb R^N).$$
\end{itemize}
\end{lemma}
\begin{proof}
To get (i) it is useful to recall the  Caffarelli-Kohn-Nirenberg inequality (see \cite{CKN}):
 {\it if $r,q\ge1$, $\frac 1r+\frac\gamma N=\frac 1q+\frac{\alpha-1}N>0 $ then for any $\phi\in C^1_c(\mathbb R^N)$
 \begin{equation}\label{ckn}
 \left\||x|^\gamma \phi\right\|_{L^r(\mathbb R^N)}\le c
 \left\||x|^{\alpha}|\nabla \phi|\right\|_{L^q(\mathbb R^N)}.
 \end{equation}
 }
 Then we apply  
\eqref{ckn} with $\gamma= 0,$  $r ={Np\over N-p}$, $q=2$ and 
$\alpha={N(2-p)\over 2p}$  
$$\(\int\limits_{\mathbb R^N}|\phi|^{Np\over N-p}\)^{(N-p)\over Np}\le c \(\int\limits_{\mathbb R^N} |x|^{N(2-p)\over p} |\nabla\phi|^2\)^{1/2}\le c \(\int_{\mathbb R^N  }|\nabla U|^{p-2}|\nabla \phi|^{2}\)^{\frac12}.$$
and the claim follows since if $p<2$ there exists a constant $c$ such that 
\begin{equation}\label{decay2}   |x|^{N(2-p)\over p } \le c {|x|^{p-2\over p-1}\over \(1+|x|^{\frac p{p-1}}\)^{N(p-2)\over p }}\ \hbox{for any}\ x\in\mathbb R.\end{equation}

To get (ii)  it is useful to recall   the weighted
 Hardy-Sobolev inequality (see for example Lemma 2.3 in \cite{DR}):\\
 {\it 
if $q\ge 1$, $s>q-N$ and $R\ge 0$ then
\begin{equation}\label{WHS}
\int_{\mathbb R^N\setminus B_R(0)}|x|^{s-q}|\varphi|^q\le c(N,q,s) \int_{\mathbb R^N \setminus B_R(0)}|x|^s |\nabla \varphi|^q\, dx\quad
 \hbox{for any}\ \varphi\in C^1_c(\mathbb R^N).\end{equation}}

Then we  apply  \eqref{WHS} with $s=-{(N-1)(p-2)\over p-1}$ and  $q=2$ (note that $s-q>-N$ because $p<N$)
$$\begin{aligned}
\int_{\mathbb R^N\setminus B_R(0)}|x|^{-{(N-1)(p-2)\over p-1}-2}|\varphi|^2&\le c  \int_{\mathbb R^N \setminus B_R(0)}|x|^{-{(N-1)(p-2)\over p-1}} |\nabla \varphi|^2 \le c(R)\int_{\mathbb R^N  }|\nabla U|^{p-2}|\nabla \phi|^{2}.\end{aligned}$$
and the claim follows since if $p>2$ for any $R>0$ there exists $c(R)$ such that 
$$  |\nabla U(x)| \ge c{|x|^\frac{N-1}{p-1}}\ \hbox{if}\ |x|\ge R.$$

 \end{proof}

Lemma \ref{mah} allows us to define the Hilbert space
$\mathcal D^{1, 2}_*(\mathbb R^N)$, which is defined as the completion of $C^1_c(\mathbb R^N)$ with respect to the norm
 $\|\phi\|:=\left( \int_{\mathbb R^N}|\nabla U|^{p-2}|\nabla\phi|^2\, dx\right)^{1/2}$ induced by the scalar product
 $$ \langle\phi,\psi\rangle:=\int_{\mathbb R^N}|\nabla U|^{p-2}(\nabla\phi,\nabla \psi)\, dx.$$
 
Now, we can look for a weak  solution $\phi\in \mathcal D^{1, 2}_*(\mathbb R^N)$ to the linear equation \eqref{linpb}, namely
\begin{equation}  \label{linpb2} \begin{aligned}&\int\limits_{\mathbb R^N}  |\nabla  U|^{p-2}(\nabla \phi,\nabla \psi), dx+(p-2)
   \int\limits_{\mathbb R^N}  |\nabla  U|^{p-4}\(\nabla  U,\nabla\phi\)\(\nabla  U,\nabla\psi\)\, dx\\
   &={Np-N+p\over N-p}  \int\limits_{\mathbb R^N}U^{Np-2N+2p\over N-p}\phi\psi \, dx\ \hbox{ for any $\psi\in \mathcal D^{1, 2}_*(\mathbb R^N).$}
\end{aligned}\end{equation} 
All the integrals involved in \eqref{linpb2} are finite. Indeed, the integrals in the L.H.S. can be easily estimated using H\"older  inequality
and Cauchy-Schwarz  inequality. 
The finiteness of the integral in the R.H.S. is more delicate and follows by the continuous embedding of the weighted space $\mathcal D^{1, 2}_*(\mathbb R^N)$ into the the weighted space 
\begin{equation}  \label{l2}\mathrm L^2_*(\mathbb R^N)=\left\{\phi\ :\ \int_{\mathbb R^N}U^{Np-2N+2p\over N-p} \phi^2\, dx<+\infty\right\},\end{equation} which is stated in the   following result.
\begin{proposition}\label{main}
There exists $C>0$ such that
\begin{equation}\label{ex1}\int_{\mathbb R^N}|\nabla U|^{p-2} |\nabla \phi|^2\, dx\ge C\int_{\mathbb R^N}U^{{Np\over N-p}-2} \phi^2\, dx\quad \hbox{for any}\ \phi\in \mathcal D^{1, 2}_*(\mathbb R^N).\end{equation}\end{proposition}

\begin{proof}
 We will prove \eqref{ex1} for any $\phi\in C^1_c(\mathbb R^N).$
 The statement will follow by a density argument. Throughout the proof $c$ will denote a  constant (possibly depending on the parameters) which may change from line to line. Although we will not estimate the constants explicitly, it will be clear from the arguments that our claims hold.\\
It is useful to remind that
$$U(x)=c {1\over \(1+|x|^{\frac p{p-1}}\)^{N(p-2)+2p\over p }}\quad \hbox{and}\quad |\nabla U(x)|^{p-2} =c{ |x|^{p-2\over p-1}\over \(1+|x|^{\frac p{p-1}}\)^{N(p-2)\over p }} .$$
We distinguish 3 cases.
 \begin{itemize}
 \item
{\em The case ${2N\over N+2}<p<2.$}
\\

We remark that since ${2N\over N+2}<p$ H\"older's inequality implies
$$ \int\limits_{\mathbb R^N}{1\over \(1+|x|^{\frac p{p-1}}\)^{N(p-2)+2p\over p }}|\phi|^2\le
\(\int\limits_{\mathbb R^N}{1\over \(1+|x|^{\frac p{p-1}}\)^{N  }}\)^{ {N(p-2)+2p \over Np}} \(\int\limits_{\mathbb R^N}|\phi|^{Np\over N-p}\)^{2(N-p)\over Np}.$$
Now, we apply Caffarelli-Kohn-Nirenberg's inequality \eqref{ckn} (with $\gamma= 0,$   $r={Np\over N-p},$ $q=2$ and 
$\alpha={N(2-p)\over 2p}$) and we get
$$\(\int\limits_{\mathbb R^N}|\phi|^{Np\over N-p}\)^{(N-p)\over Np}\le c \(\int\limits_{\mathbb R^N} |x|^{N(2-p)\over p} |\nabla\phi|^2\)^{1/2}.$$
The claim follows because of \eqref{decay2}.

\item {\em The case $1<p\le {2N\over N+2}.$}
\\

In this case ${N(p-2)+2p\over p }\le 0$ and so
$$ \(1+|x|^{\frac p{p-1}}\)^{-{N(p-2)+2p\over p }}\le c\(1+  |x|^{-{N(p-2)+2p\over (p-1) }}\).$$
Then
 $$ \int\limits_{\mathbb R^N}{1\over \(1+|x|^{\frac p{p-1}}\)^{N(p-2)+2p\over p }}|\phi|^2\le
c\int\limits_{\mathbb R^N} |\phi|^2+c\int\limits_{\mathbb R^N}|x|^{-{N(p-2)+2p\over (p-1)  }}|\phi|^2.$$
Now, we apply Caffarelli-Kohn-Nirenberg's inequality and we get (with $\gamma=0,$  $r=q=2$ and 
$\alpha=2$) 
$$\(\int\limits_{\mathbb R^N}|\phi|^{2}\)^{1/2}\le c \(\int\limits_{\mathbb R^N} |x|^{2} |\nabla\phi|^2\)^{1/2}.$$
and  also (with $\gamma=-{N(p-2)+2p\over2( p -1)},$ $r=q=2$ and 
$\alpha=1-{N(p-2)+2p\over 2(p-1)  }={-2-Np+2N\over 2(p-1) }$) 
$$\(\int\limits_{\mathbb R^N}|x|^{-{N(p-2)+2p\over (p-1)  }}|\phi|^{2}\)^{1/2}\le c \(\int\limits_{\mathbb R^N} |x|^{{-2-Np+2N\over (p-1) }} |\nabla\phi|^2\)^{1/2}.$$
Then
$$\int\limits_{\mathbb R^N}{1\over \(1+|x|^{\frac p{p-1}}\)^{N(p-2)+2p\over p }}|\phi|^2\le c\int\limits_{\mathbb R^N}
\( |x|^{2}+|x|^{{-2-Np+2N\over (p-1) }}\)
 |\nabla\phi|^2.$$
The claim follows  because if $p\le {2N\over N+2}$ it is easy to check that there exists a constant $c$ such that 
$$  |x|^{2}+|x|^{{-2-Np+2N\over (p-1) }} \le c {|x|^{p-2\over p-1}\over \(1+|x|^{\frac p{p-1}}\)^{N(p-2)\over p }}\ \hbox{for any}\ x\in\mathbb R.$$
%Indeed, if $|x|\to0$ the R.H.S. approaches $+\infty$ while the L.H.S converges to zero, since $p\le {2N\over N+2}<{2N-2\over N}.$ and
% as $|x|\to+\infty$ the R.H.S. satisfies
%$$R.H.S.\sim {|x|^{p-2\over p-1}\over \(1+|x|^{\frac p{p-1}}\)^{N(p-2)\over p }}\sim |x|^{(2-p)(N-1)\over(p-1)}\ \hbox{as} \ |x|\to+\infty,$$
%while L.H.S. is a lower order term since $1<p\le {2N\over N+2}$ implies
%$$L.H.S.\sim |x|^  {-2-Np+2N\over (p-1) }\ \hbox{and}\ {-2-Np+2N\over (p-1)}.$$

\item {\em The case $p>2.$}
\\

The proof in this case is much more delicate because the weight $|\nabla U|^{p-2}$ has different decay as $|x|\to0$ or $|x|\to\infty.$
  Let $m\ge1$ be a fixed integer.   
 We can write
$$\int_{\mathbb R^N}U^{p^*-2}\phi^2\, dx =\underbrace{\int_{\mathbb R^N\setminus B_{2^m}(0)}U^{p^*-2}\phi^2\, dx}_{(I)}+\underbrace{\int_{B_{2^m}(0)}U^{p^*-2}\phi^2\, dx.}_{(II)},$$
where $B_{2^m}(0)$ is the ball centered at the origin with radius $2^m$.\\

First, we estimate $(I)$.
We   remark that  there exists constants $c_1,\dots,c_4$ such that
\begin{equation}\label{decay}\frac{c_1}{|x|^\frac{N-p}{p-1}}\le U(x)\le \frac{c_2}{|x|^\frac{N-p}{p-1}}\ \hbox{and}\
\frac{c_3}{|x|^\frac{N-1}{p-1}}\le |\nabla U| \le \frac{c_4}{|x|^\frac{N-1}{p-1}}\ \hbox{if}\ |x|\ge 2^m.\end{equation}
Therefore
$$\begin{aligned} &\int_{\mathbb R^N\setminus B_{2^m}(0)}U^{p^*-2}\phi^2\, dx  \\ &
 \hbox{(\scriptsize we use \eqref{decay} )}\\
&\le c\int_{\mathbb R^N\setminus B_{2^m}(0)} \frac{1}{|x|^{\frac{N-p}{p-1}(p^*-2)}}\phi^2\, dx\\ &  (\hbox{\scriptsize we set $s=-\frac{N-1}{p-1}(p-2)>2-N$ and $\beta=(p^*-2)\frac{N-p}{p-1}+s-2>0$)}
\\ &=c\int_{\mathbb R^N\setminus B_{2^m}(0)} \frac{1}{|x|^{\beta}}|x|^{s-2}\phi^2\, dx\\ &
 \hbox{(\scriptsize we use \eqref{WHS} with $q=2$ and $R=2^m$ )}\\
&\le c  \int_{\mathbb R^N\setminus B_{2^m}(0)}|x|^{s}|\nabla \phi|^2\, dx\\ &
 \hbox{(\scriptsize we use \eqref{decay} )}\\ &\le c  \int_{\mathbb R^N\setminus B_{2^m}(0)}|\nabla U|^{p-2}|\nabla \phi|^2\, dx\end{aligned}$$
and so
\begin{equation}\label{(I)}
\int_{\mathbb R^N\setminus B_{2^m(0)}}U^{p^*-2}\phi^2\, dx \le c  \int_{\mathbb R^N}|\nabla U|^{p-2}|\nabla \phi|^2\, dx\end{equation}

$$ $$
Now, we estimate $(II)$. 
Firstly, given $A:=\{x\in\mathbb R^N\,\,:\,\, 1<|x|\le 2\}$, it is useful to recall the standard interpolation inequality (see for example \cite{NI})
\begin{equation}\label{step1}
\int_A \left|u-u_A\right|^m\, dx \le c \left(\int_A |\nabla u|^r\, dx\right)^{\frac m r}\end{equation} where
 $$\frac 1 m =\frac 1 r -\frac 1 N>0 \ \hbox{and}\ u_A:=\frac{1}{|A|}\int_A u\, dx.$$
Therefore, if $r=2$, $m=\frac{2N}{N-2}>0$ by  H\"older inequality  we immediately deduce
$$\int_A \left|u-u_A\right|^2\, dx \le \left(\int_A 1^{\frac{m}{m-2}}\, dx\right)^{\frac{m-2}{m}}\left(\int_A |u-u_A|^m\, dx \right)^{\frac 2 m} \le c |A|^{\frac 2 N}\int_A |\nabla u|^2\, dx.$$
Moreover, if  $\lambda>0$ and $\lambda A:=\{\lambda x\,\, :\,\, x\in A\},$
a simply  scaling gives
\begin{equation}\label{step2}
 \int_{\lambda A}|u-u_{\lambda A}|^2\, dx \le c |A|^{\frac 2 N} \lambda^{2}\int_{\lambda A}|\nabla u|^2\, dx
\end{equation}

Now, let us introduce a sequence of disjoint annuli   $A_k:=\{x\in \mathbb R^N\,:\, 2^k<|x|\le 2^{k+1}\}$ which covers the ball $B_{2^m}(0),$ namely $A_k\cap A_h=\emptyset$ for $h\neq k$ and $$B_{2^m}(0)=\bigcup_{k=-\infty}^{m-1}A_k,$$  
so that
\begin{equation}\label{5.1}
\int_{B_{2^m}(0)}U^{p^*-2}\phi^2\, dx =\sum_{k=-\infty}^m \int_{A_k}U^{p^*-2}\phi^2\, dx.
\end{equation}
We are going to estimate each term in the sum of R.H.S. of \eqref{5.1}, taking into account that
\begin{equation}\label{stimak}
  \inf_{A_k}|\nabla U|^{p-2}\ge c \frac{2^{k\frac{p-2}{p-1}}}{(1+2^{k\frac{p}{p-1}})^{\frac{N}{p}(p-2)}}.
\end{equation}
We have

\begin{equation}\label{5.2}\begin{aligned}
&\int_{A_k}U^{p^*-2}\phi^2\, dx\\ &\hbox{\scriptsize(since $p\ge2,$  $\sup\limits_{x\in\mathbb R^N}|x|^{\frac{p}{p-1}}U^{p^*-2}= L$)}\\
 &\le L \int_{A_k}|x|^{-\frac{p}{p-1}}\phi^2\, dx\\ &\le c \int_{A_k}|x|^{-\frac{p}{p-1}}|\phi-\phi_{A_k}|^2\, dx + c\int_{A_k}|x|^{-\frac{p}{p-1}}|\phi_{A_k}|^2\, dx\\ &
 \le c \int_{A_k}|x|^{-\frac{p}{p-1}}|\phi-\phi_{A_k}|^2\, dx+ c2^{k\left(-\frac{p}{p-1}+N\right)}|\phi_{A_k}|^2\end{aligned}
\end{equation}
and
\begin{equation}\label{5.3}
\begin{aligned}
\int_{A_k}|x|^{-\frac{p}{p-1}}|\phi-\phi_{A_k}|^2\, dx&\le c 2^{-k\frac{p}{p-1}}\int_{A_k}|\phi-\phi_{A_k}|^2\, dx \\
&\hbox{\scriptsize(we use \eqref{step2}   with $\lambda=2^k$)}\\ &
\le c 2^{k\left(-\frac{p}{p-1}+2\right)}\int_{A_k}|\nabla \phi|^2\, dx\\
&\hbox{\scriptsize(we use \eqref{stimak} )}\\ &
 \le c 2^{k\left(-\frac{p}{p-1}+2\right)}\frac{(1+2^{k\frac{p}{p-1}})^{\frac{N}{p}(p-2)}}{2^{k\frac{p-2}{p-1}}}\int_{A_k}|\nabla U|^{p-2}|\nabla \phi|^2\, dx\\
&= c (1+2^{k\frac{p}{p-1}})^{\frac{N}{p}(p-2)}\int_{A_k}|\nabla U|^{p-2}|\nabla \phi|^2\, dx.\\
\end{aligned}
\end{equation}
Combining \eqref{5.2} and \eqref{5.3} we get
$$
\int_{A_k}U^{p^*-2}\phi^2\, dx \le c (1+2^{k\frac{p}{p-1}})^{\frac{N}{p}(p-2)}\int_{A_k}|\nabla U|^{p-2}|\nabla \phi|^2\, dx\\+ c 2^{k\left(-\frac{p}{p-1}+N\right)}|\phi_{A_k}|^2$$
and summing upon $k$
\begin{equation}\label{5.5}\begin{aligned}
&\int_{B_{2^m}}U^{p^*-2}\phi^2\, dx\\ &=\sum_{k=-\infty}^{m-1}\int_{A_k}U^{p^*-2}\phi^2\, dx\\  & \le c \sum_{k=-\infty}^{m-1}\left[(1+2^{k\frac{p}{p-1}})^{\frac{N}{p}(p-2)}\int_{A_k}|\nabla U|^{p-2}|\nabla \phi|^2\, dx\right]\\&+ c\sum_{k=-\infty}^{m-1} 2^{k\left(-\frac{p}{p-1}+N\right)}|\phi_{A_k}|^2\\
&\le c \sum_{k=-\infty}^{m-1}\int_{A_k}|\nabla U|^{p-2}|\nabla \phi|^2\, dx+c\sum_{k=-\infty}^{m-1}\left[2^{kN\frac{p-2}{p-1}}\int_{A_k}|\nabla U|^{p-2}|\nabla \phi|^2\, dx\right]\\&+ c\sum_{k=-\infty}^{m-1} 2^{k\left(-\frac{p}{p-1}+N\right)}|\phi_{A_k}|^2\\ & \le c \int_{B_{2^m}}|\nabla U|^{p-2}|\nabla \phi|^2\, dx  + c\underbrace{\sum_{k=-\infty}^{m-1}2^{kN\frac{p-2}{p-1}}}_{ <+\infty} \int_{B_{2^m}}|\nabla U|^{p-2}|\nabla \phi|^2\, dx+\\
&+ c\sum_{k=-\infty}^{m-1} 2^{k\left(-\frac{p}{p-1}+N\right)}|\phi_{A_k}|^2\\ &
\le c  \int_{\mathbb R^N}|\nabla U|^{p-2}|\nabla \phi|^2\, dx+ c\sum_{k=-\infty}^{m-1} 2^{k\left(-\frac{p}{p-1}+N\right)}|\phi_{A_k}|^2\end{aligned}\end{equation}
It remains to estimate  the last term of \eqref{5.5}, namely
 $\sum_{k=-\infty}^{m-1} 2^{k\left(-\frac{p}{p-1}+N\right)}|\phi_{A_k}|^2.$
Now
\begin{equation*}
\begin{aligned}
& \int_{A_k\cup A_{k+1}}|\phi-\phi_{A_k\cup A_{k+1}}|^2\, dx\\ &=\frac{1}{|A_k|+|A_{k+1}|}\left(\int_{A_k}|\phi-\phi_{A_k\cup A_{k+1}}|^2\, dx+\int_{A_{k+1}}|\phi-\phi_{A_k\cup A_{k+1}}|^2\, dx\right)\\
&  \ge  \int_{A_k}|\phi-\phi_{A_k\cup A_{k+1}}|^2\, dx\\
& \ge  \left|\int_{A_k}(\phi-\phi_{A_k\cup A_{k+1}})\, dx\right|^2\\
& = \left|\int_{A_k}\phi\, dx-\frac{|A_k|}{|A_{k+1}|+|A_k|}\int_{A_k}\phi\, dx-\frac{|A_k|}{|A_{k+1}|+|A_k|}\int_{A_{k+1}}\phi\, dx\right|^2\\
&  =\frac{1}{|A_k|+|A_{k+1}| }\left||A_{k+1}|\int_{A_k}\phi\, dx-|A_k|\int_{A_k}\phi\, dx\right|^2\\
& =\frac{|A_k||A_{k+1}|}{|A_k|+|A_{k+1}| }|\phi_{A_k}-\phi_{A_{k+1}}|^2\end{aligned}
\end{equation*}
and so
\begin{equation}\label{5.7}\begin{aligned}
|\phi_{A_k}-\phi_{A_{k+1}}|^2&\le c \frac{|A_k|+|A_{k+1}|}{|A_k||A_{k+1}|} \int_{A_k\cup A_{k+1}}|\phi-\phi_{A_k\cup A_{k+1}}|^2\, dx\\
&\hbox{\scriptsize(taking into account that $|A_k|=c  2^{kN}\(2^N-1\)$ )}\\  
&\le c 2^{-kN}\int_{A_k\cup A_{k+1}}|\phi-\phi_{A_k\cup A_{k+1}}|^2\, dx\\
&\hbox{\scriptsize(we use \eqref{step2} with $\lambda=  2^{k } $ )}\\  
&\le c 2^{k(2-N)}\int_{A_k\cup A_{k+1}}|\nabla \phi|^2\, dx \\
&\hbox{\scriptsize(we use \eqref{stimak} )}\\ 
&\le c 2^{k(2-N)}(1+2^{k\frac{p}{p-1}})^\frac{N(p-2)}{p}2^{-k\frac{p-2}{p-1}}\int_{A_k\cup A_{k+1}}|\nabla U|^{p-2}|\nabla\phi|^2\, dx
\\ &\le c\(2^{k{p-Np+N\over p-1}}+ 2^{k {p-N\over p-1}}\)\int_{A_k\cup A_{k+1}}|\nabla U|^{p-2}|\nabla\phi|^2\, dx.\end{aligned}\end{equation}

Now, we use the simple fact that for any $\eta>0$ the following inequality holds
$$(a+b)^2\le (1+\eta)a^2+\(\eta+1\over\eta\)b^2\ \hbox{for any}\ a,b\in\mathbb R.$$
Then, if we    choose $\eta>0$ so that
$1+\eta=\eta_0  2^{-\frac{p}{p-1}+N}$ where $\eta_0= {\frac{2}{1+2^{-\frac{p}{p-1}+N}}}<1, $
(this is possible because
     $N-\frac{p}{p-1}>0$, since  $p\ge 2$), we get
     $$|\phi_{A_k}|^2=|\phi_{A_k}-\phi_{A_{k+1}}+\phi_{A_{k+1}}|^2\le \eta_0  2^{-\frac{p}{p-1}+N}|\phi_{A_{k+1}}|^2+ {\eta+1\over\eta}|\phi_{A_k}-\phi_{A_{k+1}}|^2.$$
 and using \eqref{5.7} we deduce
$$\begin{aligned}
&2^{k\left(-\frac{p}{p-1}+N\right)}|\phi_{A_k}|^2\\ &\le 2^{k\left(-\frac{p}{p-1}+N\right)}2^{\left(-\frac{p}{p-1}+N\right)} \eta_0|\phi_{A_{k+1}}|^2\\ &+ c 2^{k\left(-\frac{p}{p-1}+N\right)}\(2^{k{p-Np+N\over p-1}}+ 2^{k {p-N\over p-1}}\)\int_{A_k\cup A_{k+1}}|\nabla U|^{p-2}|\nabla\phi|^2\, dx\\
&=2^{(k+1)\left(-\frac{p}{p-1}+N\right)}\eta_0|\phi_{A_{k+1}}|^2 + c \(1+2^{kN{p-2\over p-1}}\)\int_{A_k\cup A_{k+1}}|\nabla U|^{p-2}|\nabla\phi|^2\, dx\end{aligned}.$$ We sum upon $k$ and we get
$$\begin{aligned}&\sum_{k=-\infty}^{m-1}2^{k\left(-\frac{p}{p-1}+N\right)}|\phi_{A_k}|^2\\ &\le \eta_0\sum_{k=-\infty}^{m-1}2^{(k+1)\left(-\frac{p}{p-1}+N\right)} |\phi_{A_{k+1}}|^2+ \sum_{k=-\infty}^{m-1}\(1+2^{kN{p-2\over p-1}}\)\int_{A_k\cup A_{k+1}}|\nabla U|^{p-2}|\nabla\phi|^2\, dx\\
&\le \eta_0\sum_{k=-\infty}^{m}2^{k\left(-\frac{p}{p-1}+N\right)} |\phi_{A_{k}}|^2+\(1+\underbrace{\sum_{k=-\infty}^{m-1}2^{kN{p-2\over p-1}}}_{<+\infty}\)\int_{\mathbb R^N}|\nabla U|^{p-2}|\nabla\phi|^2\, dx,\end{aligned}$$
which implies
\begin{equation}\label{fin1}(1-\eta_0)\sum_{k=-\infty}^{m-1}2^{k\left(-\frac{p}{p-1}+N\right)}|\phi_{A_k}|^2\le c|\phi_{A_m}|^2+c\int_{\mathbb R^N}|\nabla U|^{p-2}|\nabla \phi|^2\, dx.\end{equation}
On the other hand,  we have
\begin{equation}\label{fin2}\begin{aligned}&|\phi_{A_m}|^2\quad\hbox{\scriptsize(we use H\"older's inequality)}\\ &\le c\int_{A_m}|x|^{-\frac{(N-1)(p-2)}{p-1}}\phi^2\, dx 
\quad\hbox{\scriptsize(the annulus $A_m\subset\mathbb R^N \setminus B_{2^m(0)}$)}\\ &\le c\int_{\mathbb R^N \setminus B_{2^m(0)}}|x|^{-\frac{(N-1)(p-2)}{p-1}}\phi^2\, dx\quad\hbox{\scriptsize(we use \eqref{(I)})}\\ &\le c\int_{\mathbb R^N}|\nabla U|^{p-2}|\nabla \phi|^2\, dx.\end{aligned} \end{equation}
Finally, combining \eqref{5.5} with \eqref{fin1} (remember that $\eta_0<1$) and \eqref{fin2} 
we get
\begin{equation}\label{(II)} \int_{B_{2^m(0)}}U^{p*-2}\phi^2\, dx\le c\int_{\mathbb R^N}|\nabla U|^{p-2}|\nabla \phi|^2\, dx.  \end{equation}
 \end{itemize}
\end{proof}

\section{Proof of Theorem \ref{main1}}\label{sec-2} 

\subsection{A wave decomposition}
First of all, let us   rewrite the linear equation 
    \eqref{linpb} as
   \begin{equation}\label{linpb1}|x|^2\Delta\phi+(p-2)\sum_{i, j=1}^N\partial^2_{ij}\phi x_ix_j+\frac{(p-2)N}{1+|x|^{p\over p-1}}(\nabla\phi, x)+\gamma_{N,p}\frac{|x|^{p\over p-1}}{(1+|x|^{p\over p-1})^2}\phi=0\end{equation}
where $ \gamma_{N, p}:=N{Np-N+p \over p-1}.$
Indeed a straightforward computation shows that
    $$\begin{aligned}
    &\rm{div}\(|\nabla U|^{p-2}\nabla \phi\)+(p-2)\rm{div}\(|\nabla  U|^{p-4}\(\nabla  U,\nabla\phi\)\nabla  U\)\\
   &=|\nabla  U|^{p-2}\Delta \phi+\(\nabla |\nabla  U|^{p-2}, \nabla \phi\)\\
  &
  +(p-2) |\nabla  U|^{p-4}\(\nabla  U,\nabla\phi\)\Delta  U\\ &+(p-2)\(\nabla  U,\nabla\phi\)\(\nabla|\nabla  U|^{p-4}, \nabla  U\)\\ &+
  (p-2)|\nabla  U|^{p-4}\(\nabla\(\nabla  U,\nabla\phi\), \nabla  U\).\\
  \end{aligned}$$
and
$$\begin{aligned}
&\nabla  U=-c_{N, p}\frac{|x|^{\frac{2-p}{p-1}}x}{(1+|x|^{\frac{p}{p-1}})^{\frac N p}}\\
&|\nabla  U|^{p-4}=c_{N, p}^{p-4}\frac{|x|^{\frac{p-4}{p-1}}}{(1+|x|^{\frac{p}{p-1}})^{\frac{ N(p-4)}{ p}}}\\ 
&\left(\nabla|\nabla  U|^{p-4}, \nabla  U\right)=-c_{N, p}^{p-3}\frac{p-4}{p-1}\frac{|x|^{-\frac{2}{p-1}}}{(1+|x|^{\frac{p}{p-1}})^{\frac{ N(p-4)}{ p}+\frac N p}}+c_{N, p}^{p-3}N\frac{p-4}{p-1}\frac{|x|^{\frac{p-2}{p-1}}}{(1+|x|^{\frac{p}{p-1}})^{\frac{ N(p-4)}{ p}+\frac N p+1}}\\
&(\nabla U, \nabla\phi)=-c_{N, p}^{p-3}\frac{|x|^{\frac{2-p}{p-1}}}{(1+|x|^{\frac{p}{p-1}})^{\frac N p}}(\nabla\phi, x)\\
&\(\nabla(\nabla U, \nabla\phi), \nabla U\)=\frac{c_{N, p}^{2}}{(1+|x|^{\frac{p}{p-1}})^{\frac{ 2N}{ p}}}(\nabla\phi, x)\left[\frac{1}{p-1}|x|^{2\frac{2-p}{p-1}}-\frac{N}{p-1}\frac{|x|^{\frac{4-p}{p-1}}}{1+|x|^{\frac{p}{p-1}}}\right]\\
&\phantom{\(\nabla(\nabla U, \nabla\phi), \nabla U\)}+c_{N, p}^{2}\frac{|x|^{2\frac{2-p}{p-1}}}{(1+|x|^{\frac{p}{p-1}})^{\frac{ 2N}{ p}}}\sum_{i, j}\partial_{ij}^2\phi x_ix_j\\
&|\nabla U|^{p-4}(\nabla  U, \nabla \phi)\Delta U=\frac{c_{N, p}^{p-2}}{(1+|x|^{\frac{p}{p-1}})^{\frac{N(p-4)}{p}+\frac{ 2N}{ p}}}(\nabla\phi, x)\left[\(\frac{2-p}{p-1}+N\)|x|^{-\frac{p}{p-1}}-\frac{N}{p-1}\frac{1}{1+|x|^{\frac{p}{p-1}}}\right]\\
&(\nabla|\nabla U|^{p-2}, \nabla\phi)=\frac{c_{N, p}^{p-2}}{(1+|x|^{\frac{p}{p-1}})^{\frac{N(p-2)}{p}}}(\nabla\phi, x)\left[\frac{p-2}{p-1}|x|^{-\frac{p}{p-1}}-\frac{N(p-2)}{p-1}\frac{1}{1+|x|^{\frac{p}{p-1}}}\right]
\end{aligned}$$
where $c_{N, p}:=\alpha_{N, p}^{\frac{N-p}{p}}\frac{N-p}{p-1}.$\\

  Now, since $U$ is radial we can make a partial wave decomposition of \eqref{linpb1}, namely we can write \begin{equation}\label{phi}\phi(x)=\sum_{k=0}^\infty \phi_k(r) Y_k(\theta), \qquad \mbox{where}\quad \psi_k(r)=\int_{S^{N-1}}\phi(r, \theta)Y_k(\theta)\, d\theta,\end{equation} 
where $r=|x|$, $\theta={x\over|x|}\in S^{N-1}$ and $Y_k(\theta)$ denotes the $k$-th spherical harmonic satisfying ($\Delta_{S^{N-1}}$ stands for the Laplace-Beltrami operator)  \begin{equation}\label{harm}-\Delta_{S^{N-1}}Y_k=\lambda_k Y_k .\end{equation} It is known  that this equation has a sequence of eigenvalues 
\begin{equation}\label{lambda}
\lambda_k=k(N+k-2),\  k=0, 1, 2, \ldots,\end{equation} whose multiplicity is finite.
In particular $\lambda_0=0$ has multiplicity $1$ and $\lambda_1=N-1$ has multiplicity $N$.\\ 

Let us   write the equations satisfied by the radial functions $\psi_k.$\\
It is known that (hereafter $'$ stands for $\frac{d}{dr}$) \begin{equation}\label{lapl}\Delta(\psi_k(r) Y_k(\theta))=Y_k(\theta)\left(\psi_k''+\frac{N-1}{r}\psi_k'\right)+\frac{1}{r^2}\psi_k(r)\Delta_{S^{N-1}}Y_k(\theta).\end{equation}  Now, we have to compute the other terms in \eqref{linpb1}.
It is easy to see that
$$\partial_{x_i}\phi=\psi_k'(r)\frac{x_i}{r}Y_k(\theta)+\psi_k(r)\frac{\partial Y_k}{\partial\theta_h}\frac{\partial \theta_h}{\partial x_i}$$
and 
$$\begin{aligned}\partial^2_{x_ix_j}\phi&= \psi_k''(r)\frac{x_ix_j}{r^2}Y_k(\theta)+\psi_k'(r)\left(\frac{\delta_{ij}}{r}-\frac{x_ix_j}{r^3}\right)Y_k(\theta)+\psi_k'(r)\frac{x_i}{r}\frac{\partial Y_k}{\partial \theta_h}\frac{\partial\theta_h}{\partial x_j}\\&+\psi_k'(r)\frac{x_j}{r}\frac{\partial Y_k}{\partial\theta_h}\frac{\partial\theta_h}{\partial x_i}+\psi_k(r)\frac{\partial^2Y_k}{\partial\theta_h\partial\theta_\ell}\frac{\partial\theta_\ell}{\partial x_j}\frac{\partial \theta_h}{\partial x_i}+\psi_k(r)\frac{\partial Y_k}{\partial\theta_h}\frac{\partial^2\theta_h}{\partial x_i \partial x_j}.\end{aligned}$$
Hence 
\begin{equation}\label{grad}
(\nabla \phi, x)=\sum_{i =1}^Nx_i\partial_{x_i}\phi  =\psi_k'(r)rY_k(\theta)+\psi_k(r) \frac{\partial Y_k}{\partial\theta_h}\sum_{i=1}^N\frac{\partial \theta_h}{\partial x_i}x_i=\psi_k'(r)rY_k(\theta)
\end{equation}
and
\begin{equation}\label{dersec}\begin{aligned}
\sum_{i, j=1}^N\partial^2_{x_ix_j}\phi x_i x_j&=\psi_k''(r)r^2Y_k(\theta)+2\psi_k'(r)r\sum_{i=1}^N\frac{\partial Y_k}{\partial \theta_h}\frac{\partial\theta_h}{\partial x_i}x_i+\psi_k(r)\sum_{i, j=1}^N\frac{\partial^2Y_k}{\partial\theta_h\partial\theta_\ell}\frac{\partial\theta_\ell}{\partial x_j}x_j\frac{\partial \theta_h}{\partial x_i}x_i\\&+\psi_k(r)\sum_{i, j=1}^N\frac{\partial Y_k}{\partial\theta_h}\frac{\partial^2\theta_h}{\partial x_i \partial x_j}x_ix_j=\psi_k''(r)r^2Y_k(\theta).\end{aligned}
\end{equation}
because it holds true that
$$\displaystyle\sum_{i=1}^N\frac{\partial \theta_h}{\partial x_i}x_i=0 \ \hbox{and}\ \displaystyle\sum_{i, j=1}^N\frac{\partial^2\theta_h}{\partial x_i \partial x_j}x_ix_j=0,\ h=1, \ldots, N-1.$$

Putting together \eqref{harm}, \eqref{lapl}, \eqref{grad}  and \eqref{dersec} into \eqref{linpb1} we get the following equations for any $\psi_k,$ $k=0,1,2,\dots$,
\begin{equation}\label{linpbfinale}\psi_k''+\frac{\psi_k'}{r}\left(\frac{N-1}{p-1}+\frac{(p-2)N}{p-1}\frac{1}{1+r^{p\over p-1}}\right)-\frac{\lambda_k}{r^2}\psi_k+\gamma_{N,p}\frac{r^{\frac{p}{ p-1}-2}}{(1+r^{p\over p-1})^2}\psi_k=0,\end{equation}
which can be rewritten in a weak form as
\begin{equation}\label{linpbfinale2}\mathcal L_k(\psi_k)=0,\qquad k=0, 1, 2, \ldots,\end{equation}
where the  operator $\mathcal L_k $ is defined by
$$\begin{aligned}\mathcal L_k (\psi):&=\(r^{N-1}| U'(r)|^{p-2}\psi'\)'+ (p^*-1)r^{N-1}\(U(r)\)^{p^*-2}\psi-\lambda_kr^{N-3}| U'(r)|^{p-2}\psi.\end{aligned}
$$
Since we are concerned with solutions $\psi\in \mathcal D_*^{1,2}(\mathbb R^N)$ to the linear equation \eqref{linpb}, we will look for solutions $\psi_k$ to \eqref{linpbfinale}  or \eqref{linpbfinale2}
 in the space
$\mathcal D_k$ which is the completion of $C^1_c\([0,+\infty)\)$ with respect to the norm
$$\|\psi\|_k:=\(\int\limits_0^{+\infty} r^{N-1}| U'(r)|^{p-2}|\psi'(r)|^2\, dr+\lambda_k\int\limits_0^{+\infty} r^{N-3}| U'(r)|^{p-2}|\psi(r)|^2\, dr\)^{\frac12}.$$
 \subsection{Solving the equations $\mathcal L_k(\psi)=0$}
\begin{itemize}
\item {\it The case $k=0$.}${}$\\
We know that the function $Z_0$ defined in \eqref{Z0}  as
$$Z_0(x)=\frac{N-p}{p(p-1)}\alpha_{N, p}^{\frac{N-p}{p}}\underbrace{\frac{p-1-|x|^{\frac{p}{p-1}}}{\left(1+|x|^{\frac{p}{p-1}}\right)^{\frac N p}}}_{:=\psi_0(|x|)}$$
solves  the equation \eqref{linpb1}. 
  We claim that all the solutions in $\mathcal D_0$ to $\mathcal L_0(\psi)=0$ are given by $\psi=c\psi_0$, $c\in\mathbb R.$
 Indeed, for $k=0$ we have that $\lambda_0=0$ and a straightforward computation shows that  $\psi_0\in \mathcal D_0$ and $\mathcal L_0(\psi_0)=0.$\\
We look for a second linearly independent solution of the form $$w(r)=c(r) \psi_0(r).$$ Then we get $$c''(r)\psi_0(r)+c'(r)\left[2\psi_0'(r)+\frac{\psi_0(r)}{r}\left(\frac{N-1}{p-1}+\frac{N(p-2)}{p-1}\frac{1}{1+r^{\frac{p}{p-1}}}\right)\right]=0$$ and hence 
$$\frac{c''(r)}{c'(r)}=-2\frac{\psi_0'(r)}{\psi_0(r)}-\frac{1}{r}\left(\frac{N-1}{p-1}+\frac{N(p-2)}{p-1}\frac{1}{1+r^{\frac{p}{p-1}}}\right).$$
A direct computation shows that  $$c'(r)=A\frac{(1+r^{\frac{p}{p-1}})^{\frac{N(p-2)}{p}}}{(\psi_0(r))^2 r^{\frac{N-1+N(p-2)}{p-1}}} \ \hbox{for some}\ A\in\mathbb R\setminus\{0\}.$$
Therefore
 $$c(r)\sim B r^{\frac{N-p}{p-1}}\ \hbox{and}\ 
w (r)= c(r)\psi_0(r) \sim B  \ \hbox{as}\ r\to+\infty\ \hbox{with $B\not=0$.}$$ 
However $w\not \in \mathcal D_0$ because of Lemma \ref{mah}.\\

\item
{\it The case $k=1.$}${}$\\
We know that the function $Z_i$ defined in \eqref{Zi} as
  \begin{equation}\label{zi}\ Z_i(x)=\frac{N-p}{p-1}\alpha_{N, p}^{\frac{N-p}{p}}{x_i\over|x|}\underbrace{\frac{|x|^{\frac{1}{p-1}}}{\left(1+|x|^{\frac{p}{p-1}}\right)^{\frac N p}}}_{:=\psi_1(|x|)},\  i=1, \ldots, N\end{equation} 
  solve the equation \eqref{linpb1}. We claim that all the solutions in $\mathcal D_1$ to $\mathcal L_1(\psi)=0$ are given by $\psi=c\psi_1$, $c\in\mathbb R.$
Indeed, for $k=1$ we have that $\lambda_1=N-1$ and a straightforward computation shows that    $\psi_1\in \mathcal D_1$ and $\mathcal L_1(\psi_1)=0.$\\
As above, we look for a second linearly independent solution of the form $$w(r)=c(r) \psi_1(r).$$ Then we get $$c''(r)\psi_1(r)+c'(r)\left[2\psi_1'(r)+\frac{\psi_1(r)}{r}\left(\frac{N-1}{p-1}+\frac{N(p-2)}{p-1}\frac{1}{1+r^{\frac{p}{p-1}}}\right)\right]=0$$ 
and a direct computation shows that 
 $$c'(r)=A\frac{(1+r^{\frac{p}{p-1}})^{\frac{N(p-2)}{p}}}{(\psi_1(r))^2 r^{\frac{N-1+N(p-2)}{p-1}}}\ \hbox{for some}\ A\in\mathbb R\setminus\{0\}.$$
Therefore $$c(r)\sim B r^{\frac{N-1}{p-1}+1}\ \hbox{and}\ 
 w(r)=c(r)\psi_1(r) \sim Br\ \hbox{as}\ r\to+\infty\ \hbox{with $B\not=0$}.$$  However $w\not \in \mathcal D_1$ because of Lemma \ref{mah}.\\

 \item {\it The case $k\ge2.$}${}$\\
 We claim that all the solutions in $\mathcal D_k$ of  $\mathcal L_k (\psi)=0$ are identically zero  if $k\ge2.$
Assume there exists a function $\psi_k$ such that $\mathcal L_k (\psi_k)=0,$ i.e.  for any $r\ge0$
\begin{equation}\label{te1}
\(r^{N-1}| U'(r)|^{p-2}\psi_k'\)'+ r^{N-1}\(U(r)\)^{{Np-2N+2p\over N-p}}\psi_k-\lambda_kr^{N-3}| U'(r)|^{p-2}\psi_k=0 .\end{equation}
 We claim that $\psi_k\equiv 0$ if $k\ge2.$
We argue by contradiction. Without loss of generality, we can suppose that there exists $r_k>0$ (possibly $+\infty$) such that
$\psi_k(r)>0$ for any $r\in(0,r_k)$ and $\psi_k(r_k)=0.$ In particular, $\psi_k'(r_k)\le 0.$\\
Now, let $\psi_1(r)=U'(r)$ (see \eqref{zi})
be the solution of $\mathcal L_1(\psi_1)=0,$ i.e.  for any $r\ge0$
\begin{equation}\label{te2}
\(r^{N-1}| U'(r)|^{p-2}\psi_1'\)'+ r^{N-1}\(U(r)\)^{{Np-2N+2p\over N-p}}\psi_1-\lambda_1r^{N-3}| U'(r)|^{p-2}\psi_1=0.\end{equation}
We multiply  \eqref{te1} by $\psi_1$, \eqref{te2} by $\psi_k$, we integrate between $0$ and $r_k$, we subtract the two expressions
and we get
\begin{equation}\label{te3}
\begin{aligned}&\(\lambda_k-\lambda_1\)\int\limits_0^{r_k}r^{N-3}| U'(r)|^{p-2}\psi_k\psi_1\, dr\\ &=
\int\limits_0^{r_k}\(r^{N-1}| U'(r)|^{p-2}\psi_k'\)'\psi_1\, dr-\int\limits_0^{r_k}\(r^{N-1}| U'(r)|^{p-2}\psi_1'\)'\psi_k\, dr\\
&\hbox{\scriptsize(we integrate by part and we use that  $\psi_k(r_k)=0$)}
\\
&=r_k^{N-1}| U'(r_k)|^{p-2}\psi_k'(r_k)\psi_1(r_k)
\end{aligned}\end{equation}
and a contradiction arises when $\lambda_k>\lambda_1,$ (that is $k\ge2$), since
$\psi_k'(r_k)\le0,$ $\psi_1(r)<0$ for any $r>0$ and 
$\int\limits_0^{r_k}r^{N-3}| U'(r)|^{p-2}\psi_k\psi_1\, dr<0.$
\end{itemize}

\end{document}